\documentclass[a4paper,10pt]{article}
\title{Mixing Flows with  Homogeneous Spectrum of Multiplicity Two }
\author{V.V. Ryzhikov,  A.E. Troitskaya}
\def\eps{\varepsilon}
\date{}
\begin{document}
\def\Z{{\bf{Z}} }
\def\R{{\bf{R}} }
\maketitle
 \Large
  
\section{Introduction} 
We consider measure-preserving flows on a standard probability 
space $(X,\mu)$.   A flow $T_t$
is  mixing if for all measurable sets $A$,$B$$\subset X$
$$\mu(T_tA\cap B)\to \mu(A)\mu(B),   \ \ t\to\infty.$$

Our aim is to  obtain a mixing flow with homogeneous spectrum of multiplicity 2. This problem is connected with Rokhlin's  question
on the realization of non-simple homogeneous spectrum for an ergodic dynamical system.  Non-mixing ergodic transformations with such spectrum  appeared  in
\cite{A},\cite{R9}:  for a generic transformation $T$  the product
$T\otimes  T$  had homogeneous spectrum of multiplicity 2.
There are several generalizations of this result, see for example  
\cite{Ag}, \cite{D},\cite{DL}, \cite{KL}, \cite{KR}.
Recently the realization of  homogeneous spectrum of arbitrary multiplicity was made for mixing $\Z$-actions \cite{T}. However in case of 
mixing flows ($\R$-actions) the problem is open.

In the present note we  adapt the method of \cite{R}  to realize multiplicity 2 for  mixing   flows.
First we get non-mixing flows which could be  close to mixing ones and
possess the desired properties.  Let $\Theta$ denote the orthogonal projection 
onto the constant functions in $L_2(X,\mu)$.
\vspace{5mm}

{\bf Theorem  1.} \it  
For given $\eps$, $0<\eps<1$, let  $T_t$ be  an ergodic flow with simple spectrum such that its 
 weak closure contains operators 
$$  (1-\eps)\Theta+\frac{\eps}{a}\int_0^a{T_s}ds$$
  for all  $a>0$.
  
Then
\begin{itemize}
\item[$(1)$]
the spectral measure $\sigma$ of ${{T_t}}$ and the convolution
   $\sigma\ast\sigma$ are disjoint

\item[$(2)$] the symmetric tensor square
   $T_t\odot T_t$ has simple spectrum;
\item[$(3)$] 
   $T_t\otimes T_t$ has homogeneous spectrum of multiplicity 2.
\end{itemize}  \rm
\vskip 20pt

Then via some approximation procedure  we find a mixing  flow $T_t$ such that $T_t\odot T_t$  has  
 simple spectrum. This implies the following 
assertion.
\vspace{5mm} 

{\bf Theorem  2.} \it  
There is a mixing flow   $T_t$   such that 
   $T_t\otimes T_t$ has homogeneous spectrum of multiplicity 2.
\rm

 \section{Auxiliary  non-mixing  flows. Mixing limit flows}
 Cutting-and-staking rank-one flow construction  is determined by a parameter $h_1$,  a cut sequence $r_j$ and  a spacer sequence $\bar s_j$, 
$$ \bar s_j=(s_j(1), s_j(2),\dots, s_j(r_j-1),s_j(r_j)), $$ 
where $s_j(i)\in \R^+$.

\bf Constructions. \rm We consider a special   class $F_\eps$, $0<\eps<1$, of rank-one flows, setting ($\eps$ is fixed)  $$r_j:=j,$$

$$s_j(i):=\frac i {\sqrt{j}}$$   as $ 1\leq i \leq (1-\eps)j$,
  $$s_j(i):= \frac{i-(1-\eps)j}{\sqrt{j^3}}$$  as $ (1-\eps)j< i \leq j.$
Starting from $h_1$ let's define a sequence $h_{j}$:  $$h_{j+1}=h_jr_j +\sum_{i=1}^{r_j}s_j(i).$$
Standard "rank-one" calculations show that
$$T_{h_j}\to_w\  (1-\eps)\Theta+{\eps}I.$$
So such a  flow is $(1-\eps)$-mixing in sense of \cite{S}.
In fact the weak closure of this flow contains  all limit operators mentioned in Theorem 1.
\vspace{5mm}

\bf Lemma.   \it  For any $a> 0$ the weak closure of a  flow $T_t\in F_\eps$
contains  an operator in the form 
$$  (1-\eps)\Theta+\frac{\eps}{a}\int_0^a{T_s}ds.$$
\rm 
\vspace{3mm}

For an automorphism $T$ with simple spectrum  the existence of power limits
in the form
$$
\frac{1}{q}(I+{{T}}+{{T}}^2+\dots + {{T}}^{q-2}
+\Theta)
$$
implies  that $T\otimes T$ has homogeneous spectrum 
of multiplicity 2 \cite{R}.
 Theorem 1 is an analog of this result, its proof is similar to the proof of Theorem 4.2 from \cite{R}.

To prove Theorem 2  we consider a sequence $\eps_j\to 0$ which
tends to 0 extremely  slowly. We  consider  a rank-one flow $T_t$ setting  by a fixed $h_1$

$$r_j:=j,$$

$$s_j(i):=\frac i {\sqrt{j}}$$   as $ 1\leq i \leq (1-\eps_j)j$,
  $$s_j(i):= \frac{i-(1-\eps_j)j}{\sqrt{j^3}}$$  as $ (1-\eps_j)j< i \leq j.$

If $\eps_j\to 0$  remains constant for a very long time, then the corresponding construction   is  well approximated by flows from classes $F_{\eps}$, and  this 
   provides simple spectrum of $T_t\odot T_t$.  The latter  implies homogeneous spectrum of multiplicity 2 for  $T_t\otimes T_t$.   Because of  $\eps_j\to 0$  our flow $T_t$ becomes almost  staircase rank-one construction. In a view of \cite{Ad} (mod some modification) $T_t$ is mixing.
\large

Moscow State University

sashatro@yandex.ru,  

vryzh@mail.ru

\end{document}